\documentstyle[12pt]{article}
\newtheorem{theorem}{Theorem}\newtheorem{lemma}{Lemma}
\newtheorem{corollary}{Corollary}

\newtheorem{proposition}{Proposition}

\newtheorem{conjecture}{Conjecture}
\newtheorem{definition}{Definition}
\newtheorem{question}{Question}

\newcommand{\text}[1]{\quad\mbox{#1}\quad}

\def\beq{\begin{equation}}\def\eeq{\end{equation}}
\def\beqn{\begin{eqnarray}}\def\eeqn{\end{eqnarray}}
\def\pont{\hspace{-6pt}{\bf.\ }}

\def\qed{\ifhmode\unskip\nobreak\fi\quad\ifmmode\Box\else$\Box$\fi}

\title{Vertex covers by monochromatic pieces - a survey of results and problems}
\author{Andr\'as Gy\'arf\'as \thanks{Research supported in part by
the OTKA Grant No. K104343.}\\\\[-0.8ex]
\small   Alfr\'ed R\'enyi Institute of Mathematics \\[-0.8ex]}
\begin{document}
\maketitle
\begin{abstract}
 This survey is devoted to problems and results concerning covering the vertices of edge colored graphs or hypergraphs with monochromatic paths, cycles and other objects. It is an expanded version of the talk with the same title at the Seventh Cracow Conference on Graph Theory, held in Rytro in September 14-19, 2014.
\end{abstract}
\section{Covers by paths and cycles}

In this survey $r$-coloring always means edge-coloring with $r$ colors (traditionally red and blue when $r=2$). Some part of the material is already discussed in the 2008 survey of Kano and Li \cite{KL}.

Roughly fifty years ago, in my very first paper \cite{GGY} the following statement appeared as a footnote.

\begin{proposition}\pont \label{ggy}(Gerencs\'er, Gy\'arf\'as \cite{GGY}, 1967) The vertex set of any $2$-edge-colored complete graph $K_n$ can be partitioned into a red and a blue path.
\end{proposition}

For the sake of the rigorous, empty paths and one-vertex paths are accepted as a monochromatic path of any color in Proposition \ref{ggy}.
To prove it, suppose that $R$ and $B$ are vertex disjoint red and blue paths with endpoints $r,b$ in a complete graph with a red-blue edge coloring and $v\notin V(R)\cup V(B)$. Then either $rv$ extends $R$ to a red path or $bv$ extends $B$ to a blue path or through the edge $rb$   either $R\cup rb \cup bv, B-\{b\}$ or  $B\cup br \cup rv, R-\{r\}$ is a pair of red-blue paths covering one more vertex. Note that the path partition can be found by checking the color of at most $3(n-3)+1$ edges. We may add the two edges between the starting points and endpoints of the paths $R,B$ to get a cycle which is the union of a red and a blue path, we call it a {\em simple Hamiltonian cycle}.

\begin{corollary}\pont A $2$-colored complete graph has a simple Hamiltonian cycle.
\end{corollary}

From this one can immediately get an upper bound for the $2$-color Ramsey number of a path.
\begin{corollary}\pont  $R(P_m,P_n)\le m+n-3$ if $m,n\ge 3$.
\end{corollary}

However, the bound in the corollary is not sharp, the Ramsey number is smaller.

\begin{theorem}\pont (Gerencs\'er, Gy\'arf\'as \cite{GGY}, 1967) $R(P_m,P_n)=m+\lfloor{n\over 2}\rfloor -1$ for $m\ge n\ge 2$.
\end{theorem}

Nevertheless, in the directed case a covering result can be used to provide the Ramsey number. Let $K_n^*$ denote the complete oriented graph where each pair of vertices is joined by two edges (one in each direction). A simple Hamiltonian cycle in a $2$-colored $K_n^*$ is a directed cycle which is the union of a red and a blue path.

\begin{theorem}\pont (Raynaud \cite{RAY}, 1973) A complete oriented graph has a simple Hamiltonian cycle in every $2$-coloring of its edges.
\end{theorem}

\begin{corollary}\pont (Gy\'arf\'as, Lehel \cite{GYL}, 1973; Williamson \cite{WI}, 1973) If $m,n\ge 3$ then $R(\overrightarrow{P_m},\overrightarrow{P_n})\le m+n-3.$
\end{corollary}

Note that $R(\overrightarrow{P_m},\overrightarrow{P_n})\ge m+n-3$, shown by the following $2$-coloring of $K_{n+m-2}^*$. Color $(i,j)$ red if
$i\in \{1,2,\dots, m-1\}$ and blue if $i\in \{m,m+1,\dots m+n-2\}$. The moral is that for directed paths a cover result can provide the Ramsey number. A very recent similar example will be discussed later (Lemma \ref{proklem}).

A whole range of problems can be defined to connect Ramsey  and covering results. With S\'ark\"ozy and Selkow \cite{GYSS2011} we defined $f(n,s,r,F)$ for a family $F$  of graphs and for integers $1\le s\le r$  as the largest $m$ such that in every $r$-coloring of the edges of $K_n$ at least $m$ vertices can be covered by no more than $s$ monochromatic members of $F$. However, this survey is concentrated on the cover problem,  the smallest $s$ for which $f(n,s,r,F)=n$. Proposition \ref{ggy} implies that for the family $P$ of paths, $f(n,2,2,P)=n$.

\subsection{Cycles instead of paths - Lehel's conjecture}
Lehel  conjectured (1979) that Proposition \ref{ggy} remains true if paths are replaced by cycles. (An empty set, a vertex and an edge is accepted as a cycle.) I could prove a slightly weaker statement.

\begin{theorem}\pont(Gy\'arf\'as \cite{GY83}, 1983) The vertex set of any $2$-colored complete graph can be covered by a red and a blue cycle that intersect in at most one vertex.
\end{theorem}

It took a long time to get rid of the common vertex and prove Lehel's conjecture. First  \L uczak, R\"{o}dl and Szemer\'edi \cite{LRS} succeeded to prove the conjecture for large $n$ (with a proof using the Regularity Lemma and introducing new techniques). Then Allen \cite{A} gave a simpler proof for large $n$ without using the Regularity Lemma.

\begin{theorem}\pont (\L uczak, R\"{o}dl, Szemer\'edi \cite{LRS}, 1998; Allen \cite{A}, 2008) Lehel's conjecture is true for large enough complete graphs.
\end{theorem}

Finally Bessy and Thomass\'e found a completely elementary inductive proof that works for all graphs.

\begin{theorem}\pont \label{bt}(Bessy, Thomass\'e \cite{BT}, 2010) The vertex set of any $2$-edge-colored complete graph $K_n$ can be partitioned into a red and a blue cycle.
\end{theorem}

Recently Conlon and Stein (\cite{CS}, 2014) extended Theorem \ref{bt} (using it as a black box) to locally $2$-colored complete graphs, i.e. to complete graphs where the number of colors is arbitrary but at most two colors are used on edges incident to any vertex.

\subsection{Multicoloring}

The path (cycle) partition number of a colored graph $G$ is the minimum number of vertex disjoint monochromatic paths (cycles) whose vertices cover the vertices of $G$. These numbers are denoted by $pp(G),cp(G)$, respectively. Similarly, $pc(G),cc(G)$ denote the minimum number of monochromatic paths (cycles) whose vertices cover the vertices of $G$.

\begin{conjecture}\pont \label{pathpconj}(Gy\'arf\'as \cite{Gy}, 1989) For any $r$-colored complete graph $K$,  $pp(K)\le r$. Weaker version: $pc(K)\le r$.
\end{conjecture}

One can see easily that Conjecture \ref{pathpconj} is best possible.

\begin{proposition}\pont (Erd\H os, Gy\'arf\'as, Pyber \cite{EGYP},  1991)\label{noless} There exists $r$-colored complete graphs $K$ with $pc(K)=r$.
\end{proposition}

\noindent {\bf Proof. }  Indeed, consider a complete graph whose vertices are partitioned into sets $A_i$-s where $|A_i|=2^i-1$ for $i=1,2,\dots,r$.  The edge $xy$ gets color $$\min \{i: A_i\cap \{x,y\}\ne \emptyset\}.$$
Then color $1$ must be used to cover $A_1$ and if color $i$ is not used for some $1<i\le r$ then paths of color $1,2,\dots,i-1$ can cover at most
$$\sum_{j=1}^{i-1} 2^j=2^i-2<|A_i|$$
vertices of $A_i$. Since edges of color $j>i$ cannot cover any vertex of $A_i$, some vertex of $A_i$ remains uncovered. Thus all the $r$ colors are needed in a vertex cover by monochromatic paths. \qed

It is not obvious at all that $pc(K)$ depends on $r$ only, not on $|V(K)|$. This was established in the next result.
\begin{theorem}\pont \label{pathcovc}(Gy\'arf\'as \cite{Gy}, 1989) For any $r$-colored complete graph $K$,  $pc(K)\le r^4+r^2+1$.
\end{theorem}

Conjecture \ref{pathpconj} have been strengthened as follows.

\begin{conjecture}\pont \label{cycpconj} (Erd\H os, Gy\'arf\'as, Pyber \cite{EGYP},  1991) For any $r$-colored complete graph $K$,  $cp(K)\le r$.
\end{conjecture}

Theorem \ref{pathcovc} was extended from cover to partition, from paths to cycles with a better bound.

\begin{theorem}\pont \label{egyp}(Erd\H os, Gy\'arf\'as, Pyber \cite{EGYP}, 1991) For any $r$-colored complete graph $K$,  $cp(K)\le cr^2\log{r}$.
\end{theorem}

The proof of Theorem \ref{egyp} is based on an ``absorbing monochromatic starter'' of size $c_rn$ in an $r$-colored $K_n$ where $c_r$ is a constant depending only on $r$. It is a base cycle $C$ together with the $|C|$ different vertices, the {\em rim,} each connected to a distinct edge of $C$. Removing $C$ and its rim, one can repeatedly select monochromatic cycles in the majority color using Erd\H os - Gallai theorem, until only a small set of leftover vertices remain (relative to the size of the rim).  Finally the leftover is covered using the complete bipartite graph between the rim and the leftover part. Here P\'osa's cycle cover lemma \cite{POS} is applied: the vertex set of any graph $G$ can be covered by at most $\alpha(G)$ cycles, edges and vertices. The proof ends observing that $C$, together the unused vertices of the rim, still contains a monochromatic cycle.

Recently Conlon and Stein successfully modified the argument of the previous paragraph and obtained an extension of Theorem \ref{egyp} for local $r$-colorings, where arbitrary number of colors can be used with the condition that the edge set incident to any vertex must be colored with at most $r$ colors.

\begin{theorem}\pont \label{egyploc}(Conlon, Stein \cite{CS}, 2015) For any locally $r$-colored complete graph $K$,  $cp(K)\le cr^2\log{r}$.
\end{theorem}

Very recently Lang and Stein \cite{LS} removed the log factor from Theorem \ref{egyploc} (for large $n$).

The bound of Theorem \ref{egyp} was improved as follows, giving the best currently known bound on the cycle partition problem.

\begin{theorem}\pont (Gy\'arf\'as, Ruszink\'o, S\'ark\"ozy, Szemer\'edi \cite{GRSS}, 2006) \label{best} For any $r$-colored complete graph $K_n$ with $n\ge n(r)$,    $cp(K_n)\le 100r\log{r}$.
\end{theorem}

The gain of a factor $r$ comes from two refinements of the method outlined above. First, the Regularity Lemma is used to replace the starter with a connected structure of dense monochromatic regular pairs (lifting a dense connected monochromatic matching of the cluster graph). A monochromatic matching is {\em connected} if its edges are in the same connected component in the color of the matching. The second refinement is a more efficient covering of the leftover vertices.  The essence of a nearly optimal cover strategy is the following lemma.

\begin{lemma}\pont (Gy\'arf\'as, Ruszink\'o, S\'ark\"ozy, Szemer\'edi \cite{GRSS1}, 2006) \label{absorb} Assume $[A,B]$ is an $r$-colored complete bipartite graph, $|A|\ge r|B|$. Then $B$ can be covered by at most $r$ vertex disjoint monochromatic connected matchings.
\end{lemma}

Connected matchings are used in many Ramsey-type results, but Lemma \ref{absorb} was the first result in the context of covering by connected matchings. Further applications appear in the proofs of Theorems \ref{ass3}, \ref{BBGGS} and  \ref{BS}.

Using the bipartite $r$-color version of the Regularity Lemma, Lemma \ref{absorb} could be turned into the following result.

\begin{theorem}\pont (Gy\'arf\'as, Ruszink\'o, S\'ark\"ozy, Szemer\'edi \cite{GRSS1}, 2006) \label{main}
For every fixed $r$ there exists $n_0=n_0(r)$ such that the
following is true. Assume that the edges of a complete bipartite
graph $K(A, B)$ are colored with $r$ colors, where $|A|\ge n_0$.
If $|A| \geq 2r|B|$, then $B$ can be covered by at most $3r$
vertex disjoint monochromatic cycles.
\end{theorem}

Theorem \ref{main} shows that any covering  of all but at most ${|V(K)|\over 2r}$ exceptional vertices an $r$-colored complete graph $K$ by monochromatic cycles can be extended to a complete cover by adding at most $3r$ cycles. However, it is an open problem how to achieve a linear bound (in $r$) on the number of cycles partitioning (or even covering) the non-exceptional part.

Conjecture \ref{cycpconj} was proved asymptotically for $r=3$.

\begin{theorem}\pont (Gy\'arf\'as, S\'ark\"ozy, Szemer\'edi \cite{GYSSZ2}, 2011)\label{ass3} Apart from $o(n)$ vertices, the vertex set of any $3$-colored $K_n$ can be partitioned into $3$ monochromatic cycles.
\end{theorem}

The main idea of the proof is that for even $n$, the vertices of any $3$-colored $K_n$ can be partitioned into three monochromatic connected  matchings. Soon thereafter Pokrovskiy settled the path version, Conjecture \ref{pathpconj}, for $r=3$.

\begin{theorem}\pont (Pokrovskiy \cite{PO},  2012) For any $3$-colored complete graph $K$, $pp(K)\le 3$.
\end{theorem}

In addition, Pokrovskiy showed the following.
\begin{theorem}\pont (Pokrovskiy \cite{PO},  2012) The vertex set of any $3$-colored complete graph can be covered by $3$ (not necessarily disjoint) monochromatic paths of different colors.
\end{theorem}

It is worth noting that the two previous theorems cannot be combined. There are $3$-colorings of complete graphs where the vertex set cannot be partitioned into monochromatic paths of {\em distinct} colors. As far as I know, the first example of such colorings is discovered by K. Heinrich
\cite{H}. The simplest such $3$-coloring is obtained from the factorization of $K_4$ by replacing each vertex with a set $|A_i|=k>1$  and within each $A_i$ coloring all edges with the same color. In this coloring at most $3k$ vertices can be covered by vertex disjoint monochromatic monochromatic paths of distinct colors.

A key lemma in the proofs in \cite{PO} was the $k=1$ case of the following lemma.

\begin{lemma}\pont \label{proklem}(Pokrovskiy \cite{PO1}, 2013) For any $k$, the vertex set of any $2$-colored complete graph can be partitioned into $k$ vertex disjoint red paths and a disjoint blue complete balanced $(k+1)$-partite graph.
\end{lemma}

Lemma \ref{proklem} for general $k$ was used as a main ingredient in the proof of the following Ramsey-type result.

\begin{theorem}\pont (Pokrovskiy \cite{PO1}, 2013) For all $k$ and $n\ge k+1$, $$R(P_n,P_n^k)=(n-1)k+\left\lfloor {n\over k+1}\right\rfloor$$ where $P_n^k$ is the $k$-th power of $P_n$.
\end{theorem}

Pokrovskiy  found a counterexample to Conjecture \ref{cycpconj} for $r\ge 3$. He gave  $3$-colorings of complete graphs such that any three pairwise disjoint monochromatic cycles leave at least one vertex uncovered (but his examples can be {\em covered} by three monochromatic cycles). However, he thinks that Conjecture \ref{cycpconj} is almost true.

\begin{conjecture}\pont (Pokrovskiy \cite{PO}, 2012) In every $r$-coloring of a complete graph $K$ there are $r$ vertex disjoint monochromatic cycles covering all but $c_r$ vertices of $K$, where $c_r$ is a constant depending only on $r$.
\end{conjecture}

\subsection{Non-complete host graph}

Proposition \ref{ggy} was extended in another directions, by replacing the complete graph by sparser host graphs. We start with an old and a very recent example.

\begin{theorem}\pont(Gy\'arf\'as, Jagota, Schelp \cite{GYJS}, 1997) Assume $n\ge 5$ and $G$ is a graph obtained from $K_n$ by deleting at most $m=\lfloor {n\over 2}\rfloor$ edges. Then for every $2$-coloring of $G$, $V(G)$ can be partitioned into a red and a blue path. (Best possible, not true if $m$ is replaced by $m+1$).
\end{theorem}

\begin{theorem}\pont(Schaudt, Stein \cite{SS}, 2014) Let $G$ be a $2$-colored complete $t$-partite graph, $t\ge 3$, with no partite class of size larger than ${|V(G)|\over 2}$. Then $V(G)$ can be partitioned into a red and a blue path.
\end{theorem}

Balanced complete bipartite host graphs are also received attention.

\begin{theorem}\pont  (Gy\'arf\'as \cite{Gy}, 1989) For $r$-colored $K_{n,n}$, $pc(K_{n,n})$ is bounded by a function of $r$.
\end{theorem}

Haxell proved a stronger result, conjectured in \cite{EGYP}.

\begin{theorem}\pont (Haxell \cite{HA}, 1997) \label{hax} For any $r$-colored $K_{n,n}$, $cp(K_{n,n})$ is bounded by a function of $r$. For large $r$, $cp(K_{n,n})\le c(r\log(r))^2$.
\end{theorem}

Peng, R\"odl and Ruci\'nski \cite{PRR} mentioned that using their main result, the proof method of Theorem \ref{hax} would yield a $cr^2\log{r}$ bound. A very recent further improvement, $4r^2$ is proved by Lang and Stein \cite{LS}(it works even for local $r$-colorings).

A recent conjecture is the following.

\begin{conjecture}\pont \label{pokbipconj}(Pokrovskiy \cite{PO1}, 2014) For any $r$-colored $K_{n,n}$,   $pp(K_{n,n})\le 2r-1$ (if true best possible). This is true for $r=2$.
\end{conjecture}

It is announced in \cite{SS} that Lang, Schaudt and Stein proved Conjecture \ref{pokbipconj} for $r=3$ in asymptotic form.
The analogue conjecture is proved for infinite complete bipartite graph by D. Sokoup \cite{SOD}, see Theorem \ref{infbip} in Subsection \ref{infty}.

Host graphs with large minimum degree have been considered in the next result.

\begin{theorem} \pont \label{3/4}(Balogh, Bar\'at, Gerbner, Gy\'arf\'as, S\'ark\"ozy \cite{BBGGYS}, 2012)\label{BBGGS} For every $\eta>0$ there is $n_0(\eta)$ such that the following holds. For every $n$-vertex graph $G$ with $n\ge n_0$ and $\delta(G)>({3\over 4}+\eta)n$, in every $2$-coloring of $G$ there exists a red and a vertex disjoint blue cycle, covering all but at most $\eta n$ vertices of $G$.
\end{theorem}

A recent work of DeBiasio and Nelsen \cite{DN} introduces a new approach and strengthens Theorem \ref{3/4} by removing the exceptional set of at most $\eta n$ uncovered vertices. Letzer \cite{LET} improved that result further, proving the existence of a red-blue cycle partition in every $2$-colored $n$-vertex graph $G$ with $\delta(G)>{3n\over 4}$, provided that $n$ is large enough. 

Another result replaces the minimum degree condition of the host graph by an Ore-type condition.

\begin{theorem} \pont (Bar\'at, S\'ark\"ozy \cite{BS}, 2014) \label{BS} For every $\eta>0$ there is $n_0(\eta)$ such that the following holds. For every $n$-vertex graph $G$ with $n\ge n_0$ and such that for any two non-adjacent vertices $x,y$, $deg(x)+deg(y)\ge ({3\over 2}+\eta)n$, in every $2$-coloring of $G$ there exists a red and a vertex disjoint blue cycle, covering all but at most $\eta n$ vertices of $G$.
\end{theorem}

  Another possibility for non-complete host graph is to consider graphs with a fixed independence number. S\'ark\"ozy showed that Theorem \ref{egyp} carries over as follows.

\begin{theorem} \pont \label{alphagraph}(S\'ark\"ozy \cite{SG}, 2010) The vertex set of any $r$-colored graph $G$ with $\alpha=\alpha(G)$, $V(G)$ can be partitioned into at most $25(\alpha r)^2\log{\alpha r}$ monochromatic cycles.
\end{theorem}

Note that in case of one color ($r=1$),  a well-known result  gives the best possible vertex partition.

\begin{theorem}\pont \label{alpha}(P\'osa \cite{POS}, 1963) The vertex set of any graph $G$ can be partitioned into at most $\alpha(G)$ cycles, edges and vertices.
\end{theorem}

S\'ark\"ozy suggested that perhaps P\'osa's theorem carries over to colored graphs and strengthens Theorem \ref{alpha}.

\begin{conjecture}\pont (S\'ark\"ozy \cite{SG}, 2010) For any $r$-colored graph $G$, $cp(G)\le \alpha(G)r$.
\end{conjecture}

Although Pokrovskiy's example makes this conjecture false, it is probably true in a somewhat weaker form. At least the $2$-color case is true in the following asymptotic sense.

\begin{theorem} \pont(Balogh, Bar\'at, Gerbner, Gy\'arf\'as, S\'ark\"ozy \cite{BBGGYS}, 2012) For every $\eta>0$ and for every positive integer $\alpha$, there is $n_0(\eta,\alpha)$ such that the following holds. For every $n$-vertex graph $G$ with $n\ge n_0$ and $\alpha(G)=\alpha$, in every $2$-coloring of $G$ there exists at most $2\alpha$ vertex disjoint monochromatic cycles, covering all but at most $\eta n$ vertices of $G$.
\end{theorem}

\subsection{Hypergraphs}

A {\em loose (or linear) path} in a $k$-uniform hypergraph is a sequence of edges where only the consecutive edges meet and they meet in exactly one vertex. Here an analogue of Proposition \ref{ggy} was formulated.

\begin{conjecture}\pont(Gy\'arf\'as, S\'ark\"ozy \cite{GyS}, 2013) In any $2$-colored complete $k$-uniform hypergraph there exists a vertex disjoint red and a blue loose path  covering all but at most $k-2$ vertices. (There is an example showing that one cannot avoid $k-2$ uncovered vertices.)
\end{conjecture}

A weaker form of this conjecture, with $2k-5$ uncovered vertices instead of $k-2$, is proved in \cite{GyS}, thus the case $k=3$ is settled.

A {\em tight path} in a $k$-uniform hypergraph is a sequence of at least $k$ vertices and with all edges that have $k$ vertices in consecutive positions in the sequence. It is possible that Proposition \ref{ggy} extends to tight paths, but we could not even decide the first step.

\begin{question}\pont(Gy\'arf\'as, S\'ark\"ozy \cite{GyS}, 2013)
Is it possible to partition the vertex set of every $2$-colored complete $3$-uniform hypergraph into a red and a blue tight path?
\end{question}

However, in the case of Berge paths, when only distinct vertices and distinct edges (covering consecutive pairs of vertices) are required in the definition of paths, Proposition \ref{ggy} can be extended to any number of colors.

\begin{theorem}\pont(Gy\'arf\'as, S\'ark\"ozy \cite{GyS}, 2013)  The vertex set of any $r$-uniform $r$-colored complete hypergraph $K$ can be partitioned into monochromatic Berge paths of different colors. Consequently, at most $r$ Berge paths cover $V(K)$.
\end{theorem}

The proof method in \cite{EGYP} was extended to bound the cycle partition number for loose cycles of hypergraphs.

\begin{theorem}\pont\label{loosepart}(Gy\'arf\'as, S\'ark\"ozy \cite{GyS}, 2013) There exists $c(r,k)$ such that the vertex set of any $r$-colored complete $k$-uniform hypergraph can be covered by at most $c(r,k)$ vertex disjoint monochromatic loose cycles.
\end{theorem}

S\'ark\"ozy gave a good bound on $c(r,k)$ from the strong hypergraph Regularity Lemma of R\"odl and Schacht \cite{RS}, combined with the hypergraph Blow-up lemma of Keevash \cite{KE}.

\begin{theorem}\pont(S\'ark\"ozy \cite{SA1}, 2014) For $c(r,k)$ in Theorem \ref{loosepart}, we have $c(r,k)\le 50rk\log{rk}$.
\end{theorem}

The most general existence theorem, combining Theorems \ref{alphagraph} and \ref{loosepart} is obtained recently.

\begin{theorem} \pont(Gy\'arf\'as, S\'ark\"ozy \cite{GYS14}, 2014)\label{exists} For every $r$-coloring of a $k$-uniform hypergraph $H$ with $\alpha=\alpha(H)$, $V(H)$ can be partitioned into at most $c(r,k,\alpha)$ monochromatic loose cycles.
\end{theorem}

Concerning sharp results, even for one color, the following possible extension of P\'osa's theorem is suggested.

\begin{conjecture} \pont(Gy\'arf\'as, S\'ark\"ozy \cite{GYS14}, 2014) Every $k$-uniform hypergraph can be partitioned into at most $\alpha(H)$ loose cycles and parts of hyperedges. (The case $k=2$ is P\'osa's theorem.)
\end{conjecture}

We showed that the conjecture is true if loose cycles replaced by weak cycles (sequence of edges in which only the cyclically consecutive edges intersect) or by loose paths. It seems reasonable that Theorem \ref{exists} can be extended to tight cycles (where all cyclically consecutive $k$-sets form the edges).

\begin{conjecture}\pont There exists $c(r,k)$ bounding the tight cycle partition number of any $r$-colored complete $k$-uniform hypergraph.
\end{conjecture}

\subsection{Infinite versions}\label{infty}

For infinite complete graphs I think that the following result is ``perfect''.

\begin{theorem}\pont (Rado \cite{RA}, 1987) Countably infinite $r$-colored complete graphs can be partitioned into finite or one-way infinite monochromatic paths of different colors. Consequently, at most $r$ monochromatic paths cover the vertex set.
\end{theorem}

As in the finite case, $r-1$ monochromatic paths cannot always cover, the coloring of Proposition \ref{noless} with $|A_r|=\infty$ shows that.  Rado's proof can be easily extended to hypergraphs.

\begin{theorem}\pont\label{infloose}(Gy\'arf\'as, S\'ark\"ozy \cite{GyS}, 2013) The vertex set of any countably infinite $k$-uniform $r$-colored hypergraph can be partitioned into monochromatic loose paths of different colors.
\end{theorem}

The following result (in fact it is a special case of a more general result) extends Theorem \ref{infloose} to tight paths.

\begin{theorem}\pont (Elekes, Soukup, Soukup, Szentmikl\'ossy \cite{ESSSZ}, 2014) The vertex set of any countably infinite $k$-uniform $r$-colored hypergraph can be partitioned into monochromatic tight paths of different colors.
\end{theorem}

The infinite version of Conjecture \ref{pokbipconj}  was proved by D. T. Soukup.

\begin{theorem}\pont (Soukup \cite{SOD}, 2014)\label{infbip} For any $r$-colored balanced countably infinite complete bipartite graph $B$,  $pp(B)\le 2r-1$.
\end{theorem}

\subsection{A misunderstanding - when red and blue cannot mix}
I mentioned once to Paul Erd\H os that for any $2$-coloring of $K_n$, $pc(K_n)\le 2$ (Proposition \ref{ggy}) and he said he did not believe it. It soon turned out that Paul thought that the covering paths must be monochromatic {\em in the same} color!   Indeed, with this definition the path cover number depends on $n$.

\begin{theorem}\pont (Erd\H os, Gy\'arf\'as \cite{EGy}, 1995) The vertex set of any $2$-colored $K_n$ can be covered by at most $2\sqrt{n}$ monochromatic paths of the same color.
\end{theorem}

\begin{conjecture}\pont (Erd\H os, Gy\'arf\'as \cite{EGy}, 1995) The vertex set of any $2$-colored $K_n$ can be covered by at most $\sqrt{n}$ monochromatic paths of the same color. This would be best possible.
\end{conjecture}

\section{Covers by other monochromatic subgraphs}

\subsection{Covers with graphs of bounded degree}

The cycle cover problem can be extended to coverings by regular graphs.

\begin{theorem} \pont(S\'ark\"ozy, Selkow, Song \cite{SAS1},  2013) Every $r$-colored complete graph $K_n$ can be covered by at most $100r\log{r}+2rd$ vertex disjoint connected monochromatic $d$-regular graphs and vertices for $n\ge n_0(r,d)$.
\end{theorem}

The most general result into this direction is found by Grinshpun and S\'ark\"ozy.  They call a sequence of graphs $F=\{F_1,F_2,\dots\}$  $\Delta$-bounded if each $F_i$ has $i$ vertices and maximum degree at most $\Delta$.

\begin{theorem} \pont (Grinshpun, S\'ark\"ozy \cite{GSA}, 2014) There is an absolute constant $C$ such that for every $\Delta$ and every $\Delta$-bounded sequence $F$, any $2$-colored complete graph can be partitioned into at most $2^{C\Delta \log{\Delta}}$ vertex disjoint monochromatic members of $F$. If all $F_i$s are bipartite graphs, $2^{C\Delta}$ members suffice and this is best possible apart from the constant in the exponent.
\end{theorem}

\begin{corollary} \pont \label{kpower}(Grinshpun, S\'ark\"ozy \cite{GSA}, 2014) The vertex set of every $2$-colored complete graph can be covered by at most $2^{ck\log{k}}$ vertex disjoint monochromatic $k$-th power of a cycle.
\end{corollary}

The authors conjecture sharper bounds for the covering number with powers of cycles.

\begin{conjecture}\pont (Grinshpun, S\'ark\"ozy \cite{GSA}, 2014) The vertex set of every $2$-colored complete graph can be covered by at most $ck$ vertex disjoint monochromatic $k$-th power of a cycle.
\end{conjecture}

\subsection{Covers with connected pieces - Ryser's conjecture}

A subset $S$ of vertices in a colored hypergraph is called {\em connected  in color $i$} if the edges of color $i$ within $S$ form a connected subhypergraph. $S$ is {\em connected} if it is connected in some color.

\begin{definition}\pont The connected cover (partition) number of a colored graph or hypergraph $H$ is the minimum number of (vertex disjoint) connected subsets covering $V(H)$. These numbers are denoted by $conc(H),conp(H)$, respectively.
\end{definition}

A conjecture of Ryser (1971) that appeared in \cite{HE} was formulated in \cite{GYL2} as follows.

\begin{conjecture} \pont \label{rysconj} For every $r$-colored graph $G$, $conc(G)\le (r-1)\alpha(G)$ ($r\ge 2$).
\end{conjecture}
Ryser's conjecture is open for $\alpha(G)=1$ except for $r\le 5$. The case $r=2$ is equivalent with K\"onig's theorem and for  $r=3$ it was proved by Aharoni (\cite{AH}, 2001). Ryser's conjecture is sharp if $r-1$ is a prime power and $|V(G)|$ is divisible by $(r-1)^2$.

\subsection{Covers with connected pieces - complete bipartite host}

Ryser's conjecture for $\alpha(G)=1$ states that $r$-colored complete graphs can be covered by at most $r-1$ monochromatic connected subsets. 
For complete bipartite graphs an analogue conjecture is the following.

\begin{conjecture}\pont \label{rysbip}(Gy\'arf\'as, Lehel \cite{GYL2}, 1977) For every $r$-colored complete bipartite graph $G$, $conc(G)\le 2r-2$.
\end{conjecture}

A non-trivial example (\cite{GYL2}, \cite{CFGYLT}) shows that $2r-2$ is best possible in Conjecture \ref{rysbip}.  It is worth comparing Conjecture \ref{rysbip} with Conjecture \ref{pokbipconj} and Lemma \ref{absorb}.
 
Conjecture \ref{rysbip} can be reduced to $r$-colorings where each color class is the union of vertex disjoint complete bipartite graphs, covering all vertices. Such reduction is not known for Ryser's conjecture.

\begin{theorem} \pont (Chen, Fujita, Gy\'arf\'as, Lehel, T\'oth \cite{CFGYLT}, 2013) Conjecture \ref{rysbip} is true for $r\le 5$.
\end{theorem}

\subsection{Covers with connected pieces - hypergraphs}
Although Ryser's conjecture is open for complete graphs ($\alpha(G)=1, r\ge 6$), its analogue for complete $k$-uniform hypergraphs is solved by Zolt\'an Kir\'aly with a short elegant proof.

\begin{theorem} \pont (Kir\'aly \cite{KI}, 2010) If $k\ge 3$, for every $r$-colored complete $k$-uniform hypergraph $H$, $conc(H)\le \lceil{r\over k}\rceil$. The bound is sharp for $k\ge 3, r\ge 1$.
\end{theorem}

For non-complete $k$-uniform hypergraphs $H$ ($\alpha(H)\ge k$) only the initial steps are known and there is no plausible conjecture for the best estimate of $conc(H)$. Here $\alpha(H)$ is the cardinality of the largest vertex set that does not contain completely any edge of $H$.

\begin{theorem}\pont (Fujita, Furuya, Gy\'arf\'as, T\'oth \cite{FFGYT2}, 2014) For $k\ge 2$ every $k$-colored $k$-uniform hypergraph $H$ with $\alpha(H)=k$ satisfies $conc(H)\le 2$. For $k\ge 3$ every $(k+1)$-colored $k$-uniform hypergraph $H$ with $\alpha(H)=k$ satisfies $conc(H)\le 3$. Both results are sharp.
\end{theorem}

\subsection{Partitions into connected pieces}

It is possible that Conjecture \ref{rysconj} can be strengthened from cover to partition. For the case $\alpha(G)=1$, i.e. for complete graphs, this was formulated in \cite{EGYP}.

\begin{conjecture}\pont(Erd\H{o}s,  Gy\'{a}rf\'{a}s, Pyber \cite{EGYP}, 1991) For every $r$-colored complete graph $K$, $conp(K)\le r-1$.
\end{conjecture}

For $r\le 3$ this conjecture was proved in \cite{EGYP}. The next result ''almost proves'' it.

\begin{theorem} \pont(Haxell, Kohayakawa \cite{HAKO}, 1996) If $n$ is large compared to $r$ then for every $r$-colored complete graph $K$, $conp(K)\le r$.
\end{theorem}

For general $\alpha$ and for $k$-uniform hypergraphs, the following bound is known.

\begin{theorem}\pont (Fujita, Furuya, Gy\'arf\'as, T\'oth \cite{FFGYT1}, 2012) For every $r$-colored non-trivial $k$-uniform hypergraph $H$, $conp(H)\le \alpha(H)-r+2$.
\end{theorem}

When $r=2$ this result gives the following (an extension of K\"onig's theorem).

\begin{corollary}\pont (Fujita, Furuya, Gy\'arf\'as, T\'oth \cite{FFGYT1}, 2012) For every $2$-colored non-trivial $k$-uniform hypergraph $H$, $conp(H)\le \alpha(H)$.
\end{corollary}

\noindent {\bf Acknowledgement.} Thanks for G\'abor S\'ark\"ozy for his remarks on the manuscript. The author is also indebted to referees for valuable remarks.

\end{document}